\documentclass[11pt]{article}
\usepackage{graphicx}
\usepackage{amsmath,amssymb,amsthm,amsfonts}
\usepackage{hyperref}
\hypersetup{hidelinks,
	colorlinks=true,
	filecolor=red,      
	urlcolor=green,
	citecolor=blue,
	pdfstartview=Fit,
	breaklinks=true}

\newtheorem{thm}{Theorem}[section]

\newtheorem{lem}[thm]{Lemma}

\theoremstyle{definition}

\theoremstyle{remark}

\numberwithin{equation}{section}

\DeclareMathSymbol{\C}{\mathalpha}{AMSb}{"43}

\newcommand{\R}{\mathbb{R}}

\begin{document}
\title{Existence and  Asymptotic Behavior of Minimizers for Rotating Bose-Einstein Condensations  in Bounded Domains}
\author{
	Yongshuai Gao\thanks{School of Mathematics and Statistics, Central China Normal University, 
		Wuhan 430079, People's Republic of China.  Email: \texttt{ysgao@mails.ccnu.edu.cn}. }, 
	\, Shuai Li\thanks{College of Informatics, Huazhong Agricultural University,  
		Wuhan 430070, People's Republic of China.  Email: \texttt{lishuai@mail.hzau.edu.cn}. }
		% S.  Li is partially supported by NSFC grant 12371113 and 11901223.
	\, and
	\, Peiye Zhong\thanks{School of Mathematics and Statistics, Central China Normal University,  
		Wuhan 430079, People's Republic of China.  Email: \texttt{peiyezhong@mails.ccnu.edu.cn}. }
}

\date{\today}

\smallbreak \maketitle

\begin{abstract}
This paper is concerned with the existence and mass concentration behavior of minimizers for rotating Bose-Einstein condensations (BECs) with  attractive interactions in a bounded domain $\mathcal{D}\subset \R^2$. 
It is shown that,  there exists a finite constant  $a^*$,  denoting mainly the critical number of  bosons  in the system, 
  such that  
the least energy $e(a)$ admits minimizers if and only if  $0<a<a^*$, 
no matter  the trapping potential $V(x)$ rotates at any velocity  $\Omega\geq0$. 
This is quite different from the rotating BECs in the whole plane case, where the existence conclusions depend on the value of  $\Omega$ (cf. \cite[Theorem 1.1]{GLY}). 
Moreover, by establishing the refined  estimates of the  rotation term and the least energy,  we  also analyze the mass concentration behavior of minimizers in a harmonic potential as $a\nearrow a^*$.  
\end{abstract}


\vskip 0.2truein

\section{Introduction}
Bose-Einstein condensation (BEC) has been investigated widely since it was first realized in 1990s \cite{AEMWC,Ketterle}. 
In physical experiments of rotating BECs, a large number of  (bosonic) atoms confined in the rotating traps simultaneously occupy the ground state of the system as temperature is below a critical value (cf. \cite{ARVK,A,CR,DGPS,Fetter}). 
These rotating BECs display various interesting quantum phenomena, such as the appearance of quantum vortices \cite{ARVK,A,CC,Fetter} and the center-of-mass rotation \cite{ARVK,Fetter,LCS}. 
Over the last two decades,  the study of  rotating BECs is always a core topic of physicists and mathematicians both domestically and internationally, see  \cite{CC,CD,CPRY,CRY,GLP,GLY,HMDBB,LNR} and the references therein.

In this paper, we are interested in the ground states of rotating BECs
with attractive interactions in a bounded domain $\mathcal{D}\subset\R^2$, which can be described by minimizers of the following complex-valued
variational problem (cf. \cite{ANS,BC,Fetter,LNR,LS}): 
\begin{equation}\label{eq:1.1}
	e(a):=\inf_{\substack{\{u\in H^1_0(\mathcal{D},\C),\|u\|_2^2=1\}}} E_a(u),
\end{equation}
where the Gross-Pitaevskii (GP) functional $E_a(u)$ is given by
\begin{equation}\label{eq:1.2}
	E_a(u):=\int_{\mathcal{D}} \big(|\nabla u|^2+V(x)|u|^2\big)dx-\frac{a}{2}\int_{\mathcal{D}} |u|^4dx-\Omega \int_{\mathcal{D}} x^\bot\cdot\big(iu,\nabla u \big)dx.
\end{equation}
Here $x=(x_1,x_2)\in\mathcal{D}\subset\R^2$, $x^\bot=( -x_2,x_1)$, $(iu,\nabla u):=i(u\nabla \bar u-\bar u\nabla u)/2$, and $\mathcal{D}$ is a bounded domain satisfying the interior ball condition and $\partial\mathcal{D}\in C^1$. The parameter $a>0$ in \eqref{eq:1.2} characterizes the absolute product of the scattering length $\nu$ of the two-body interaction times the number $N$ of particles in the condensates, and $\Omega\geq0$ describes the rotational velocity of the trap $V(x)$.  

Recall from \cite[Theorem 7.21]{LL} (see also  (2.3) in \cite{EL1989}) that, for any $u\in H^1(\R^2,\C)$ and $\mathcal{A}\in L^2_{loc} (\R^2,\R^2)$, the following diamagnetic inequality holds: 
\begin{equation}\label{eq:1.7}
	|(\nabla-i\mathcal{A})u|^2
	=|\nabla u|^2-2\mathcal{A}\cdot \big(iu,\nabla u \big)+\mathcal{A}^2|u|^2
	\geq |\nabla|u||^2 \  a.e. \  \mbox{in} \  \R^2.
\end{equation} 
Set  $\mathcal{A}=\frac{\Omega}{2}x^\bot$, 
and  the GP functional $E_a(u)$ can also be rewritten as
\begin{equation}\label{exp:Eau.2}
	E_a(u)=\int_{\mathcal{D}}\Big |\Big(\nabla-i \frac{\Omega}{2}x^{\bot}\Big)u\Big|^2dx
	+\int_{\mathcal{D}}V_\Omega(x) |u|^2dx-\frac{a}{2}\int_{\mathcal{D}} |u|^4dx, 
\end{equation}
where $V_\Omega(x)$ is defined by 
\begin{equation}\label{def:V_Omega}
	V_\Omega(x):=V(x)-\frac{\Omega^2}{4}|x|^2, \,\,\, x\in\mathcal{D}. 
\end{equation}

When there is no rotation for $V(x)$, i.e. $\Omega=0$, one can note from the diamagnetic inequality \eqref{eq:1.7} that
\begin{equation}\label{ineq: dia.pri}
	\text{$|\nabla u|\geq|\nabla|u||$ $a.e.$ in $\R^2$, }
\end{equation}
which indicates that all  minimizers of \eqref{eq:1.1} are essentially real-valued. 
In the last decades,   existence, refined limit  behavior,   symmetry breaking and   local uniqueness of  minimizers for  problem \eqref{eq:1.1} with $\Omega =0$ in the whole space $\R^2$ was widely studied, see  \cite{BC,GLW,GS,GWZZ,GZZ,Zhang} and the references therein.  
Denote by $w(x)$  the unique (up to translations) positive radially symmetric solution of the following nonlinear scalar field equation (cf.  \cite{GNN,K,W1983})
\begin{equation}\label{eq:1.3}
	-\Delta u+u-u^3=0 \,\ \text{in} \,\ \R^2, \,\ \text{where} \,\ u\in H^1(\R^2,\R). 
\end{equation}  
It was proved in \cite{BC, GS} that,  
$e(a)$ in $\R^2$ admits minimizers if and only if  $0\leq a<a^*:=\|w\|_2^2$. 
Moreover, by establishing corresponding refined energy estimates,  concentration behavior and local uniqueness of minimizers as $a\nearrow a^*$ were also analyzed by  \cite{GLW,GS,GWZZ,GZZ} under different types of trapping potentials. 
These above mentioned results show that the 
analysis of problem \eqref{eq:1.1} in $\R^2$ makes full use of the following classical  Gagliardo-Nirenberg (GN) inequality in $\R^2$ (cf. \cite[Theorem B]{W1983}):  
\begin{equation}\label{ineq:GN.R2}
	\int_{\R^2}|u|^4dx\leq 
	\frac{2}{a^*}\int_{\R^2}|\nabla u|^2dx\int_{\R^2}|u|^2dx,
	\ \ u\in H^1(\R^2,\R),
\end{equation}
where $\frac{2}{a^*}$ is the best constant,  and up to a scaling, the identity is achieved at  $u=w(x)$.

More recently, problem \eqref{eq:1.1} with $\Omega=0$ in  a bounded domain was studied in \cite{LZ}.  In order to analyze the existence and asymptotic behaviors of the minimizers,  the authors improve \eqref{ineq:GN.R2} into a  new type of  GN inequality:   
\begin{equation}\label{eq:1.6}
	\int_{\mathcal{D}} |u|^4dx\leq \frac{2}{a^*} \int_{\mathcal{D}} |\nabla u|^2dx\int_{\mathcal{D}} |u|^2dx,\  \     u\in H_0^1(\mathcal{D},\R), 
\end{equation}
where  the best constant $\frac{2}{a^*}$ cannot be attained.  

When the trapping potential $V(x)$ rotates at a velocity  $\Omega>0$,  minimizers of  \eqref{eq:1.1} are complex-valued. Recently, some qualitative properties for minimizers of  \eqref{eq:1.1}  in $\R^2$ have  attracted the attention of scholars gradually,  such as the existence and nonexistence, asymptotic behaviors, local uniqueness,  stability and free-vortex and so on,  see  \cite{ANS,BC,Guo,GLP,GLY,LNR} and the references therein. 
In detail, applying  the diamagnetic inequality \eqref{eq:1.7} and the GN inequality \eqref{ineq:GN.R2}, 
it was proved in \cite{GLY,LNR}  that  \eqref{eq:1.1} in $\R^2$ admits minimizers if and only if $0<a<a^*$, provided that the rotating velocity $\Omega$ is smaller than a critical value $\Omega^*$. 
Furthermore, using the method of inductive symmetry,  the uniqueness and free-vortex of minimizers for \eqref{eq:1.1} in $\R^2$  as $a\nearrow a^*$ were proved in \cite{GLY} for some suitable class of radial potential $V(x)=V(|x|)$. 
Subsequently,  more general results on the local uniqueness and nonexistence of vortices
for minimizers in some non-radially symmetric trap were also analyzed by \cite{Guo} and \cite{GLP}, respectively.

Motivated by the works mentioned above, the aim of this paper is to 
 investigate the existence and nonexistence, and asymptotic behavior of minimizers for  \eqref{eq:1.1} in a bounded domain $\mathcal{D}$.  From the physical point of view, we always assume that the trapping potential  $V(x)$ is harmonic, $i.e.$, 
\begin{equation}\label{eq:1.8}
	V(x)=x_1^2+\Lambda x_2^2, 
\end{equation}
where $\Lambda\in\R^+$ and $x=(x_1,x_2)\in \mathcal{D}$. 
Moreover,  set 
\begin{equation}\label{cond:V_Omega}
	Z:=\big\{x\in\bar{\mathcal{D}}:  V_\Omega(x)=\inf\limits_{y\in\bar{\mathcal{D}}}V_\Omega(y)\big\},  
\end{equation}
where $\bar{\mathcal{D}}:=\mathcal{D}\cup\partial\mathcal{D}$, 
and $V_\Omega(x)$ on $\partial\mathcal{D}$ is defined by 
\[
V_\Omega(x)\big|_{x=(x_1,x_2)\in\partial\mathcal{D}}:=x_1^2+\Lambda x_2^2-\frac{\Omega^2}{4}|x|^2. 
\]
One can observe that $Z$ denotes the set of  minimum points of $V_\Omega(x)$ in $\bar{\mathcal{D}}$.

 Under the above assumptions, we shall now  state our first result on the  existence and nonexistence of minimizers for problem \eqref{eq:1.1} as the following theorem. 

\begin{thm}\label{thm.exist}
 Suppose $V(x)$ satisfies \eqref{eq:1.8} and $w(x)$ denotes the unique positive radially symmetric solution of \eqref{eq:1.3}. Set $a^*:=\|w\|_2^2$. Then 
\begin{enumerate}
\item  If $0\leq\Omega<\infty $ and  $0< a<a^*$, there exists at least one minimizer for e(a).
\item  If $0\leq\Omega<\infty $ and  $a\geq a^\ast$, there is no minimizer for e(a).
\end{enumerate}
\end{thm}

Theorem \ref{thm.exist} indicates that, for any $\Omega\in[0,\infty)$, $e(a)$ admits minimizers if and only if $0< a<a^*$. 
This  is quite different from the existence results for  $e(a)$ in the whole space $\R^2$. 
As stated before (see also  \cite[Theorem 1.1]{GLY}), there exist a critical velocity $0<\Omega^*\leq\infty$ such that  $e(a)$ in $\R^2$ admits minimizers if and only if $0< a<a^*$,   provided that  $\Omega<\Omega^*$.  
The reason for this different conclusion here is the boundness of the domain $\mathcal{D}$, which leads to that the rotational term $-\Omega\int_{\mathcal{D}} x^\bot\cdot(iu,\nabla u)dx$ is bounded from below for any given $a\in(0,a^*)$.  
Besides, note from the proof of Theorem \ref{thm.exist} in Section 2, and we remark that  Theorem \ref{thm.exist} also holds for some  more general case where the trapping potential $V(x)\in L^\infty(\mathcal{D})$ is  bounded from below.

Suppose $u_a\in H_0^1(\mathcal{D},\C)$ is a minimizer of $e(a)$ for $0< a <a^\ast$, and then by the variational argument, $u_a$  satisfies the following Euler-Lagrange equation
 \begin{equation}\label{eq:3.4}
 	\begin{cases}
 		-\Delta\,u_a+V(x)u_a+i\Omega x^\bot \cdot \nabla u_a=\mu_a u_a+a|u_a|^2 u_a,  &\text{ in $\mathcal{D}$ },\\
 		u_a=0,                    & \text{ on $\partial\mathcal{D}$, } 
 	\end{cases}
 \end{equation}
 where $\mu_a\in\R$ is the associated Lagrange multiplier and
 \begin{equation}\label{eq:3.5}
 	\mu_a=e(a)-\frac{a}{2}\int_{\mathcal{D}} |u_a|^4dx,\    0<a<a^\ast.
 \end{equation}
Our second result is concerned with the mass concentration behavior of minimizers for problem \eqref{eq:1.1}, where the trap potential $V(x)$ satisfies \eqref{eq:1.8} with $\Lambda>1$ and $\Omega\leq2$.  
One can check from \eqref{def:V_Omega} and \eqref {eq:1.8} that,   $V_{\Omega}(x)\geq0$ and the origin is the unique  minimum point of $V_{\Omega}(x)$  if  $\Lambda>1$ and $\Omega\leq2$.    Hence, we shall always suppose $0\in\mathcal{D}$ as we analyze the limit behavior of minimizers in the following theorem. 

\begin{thm}\label{thm.1.2}
Suppose the origin is a inner point of $\mathcal{D}$, and 
$V(x)$ satisfies \eqref{eq:1.8} with $\Lambda>1$ and $\Omega\leq2$.  
Let $u_a$ be a complex-valued minimizer of $e(a)$,  and then we have  several conclusions as follows.  
\begin{enumerate}
	\item  $|u_a|$ admits a unique local maximum point $x_a$ as $a$ is close enough to $a^*$,  and $x_a$ satisfies  
	\begin{equation}\label{lim.xa}
		\frac{ |x_a|}{(a^\ast-a)^{\frac{1}{4}}}\to0  \text{ as } a  \nearrow a^\ast.
	\end{equation}
	\item  Set $\tilde{u}_a:=\sqrt{a^\ast}\frac{(a^\ast-a)^{\frac{1}{4}}}{\lambda} u_a\Big(\frac{(a^\ast-a)^{\frac{1}{4}}}{\lambda} x+x_a\Big)e^{-i(\frac{\Omega}{2\lambda}(a^\ast-a)^{\frac{1}{4}} x \cdot x_a^\bot-\theta_a)}$, where  $\theta_a\in[0,2\pi)$ is a proper constant and 
 \begin{equation}\label{eq:1.10}
	\lambda:=\Big[\int_{\R^2} V(x)w^2dx\Big]^\frac{1}{4}>0.
\end{equation}
Then $\tilde{u}_a$ satisfies 
\begin{equation}\label{eq:1.9}
	\tilde{u}_a\to w(x)\  \mbox{  uniformly in $L^\infty(\R^2,\C)$ as $a  \nearrow  a^\ast$},  
\end{equation}
where $w$ is the unique positive radially symmetric solution of \eqref{eq:1.3}. 
%	\item $e(a)$ satisfies 
%\begin{equation}\label{lim:ea}
%	\frac{e(a)}{(a^\ast-a)^{\frac{1}{2}}}\to\frac{2\lambda^2}{a^\ast}  \text{ as } a  \nearrow a^\ast. 
%\end{equation}
\end{enumerate}
\end{thm}

Note from the proof of Theorem \ref{thm.1.2} in Section \ref{sec.MC}  that the limit \eqref{eq:1.9} holds in $L^\infty(\R^2,\C)$ in the sense that
$u_a(x)\equiv0$ for all $x\in\mathcal{D}^c:=\R^2\setminus\mathcal{D}$. 
Theorem \ref{thm.1.2} shows that  minimizers of problem \eqref{eq:1.1} must concentrate at the origin as $a\nearrow a^*$. The rates of the unique maximum point of $|u_a|$ tending to the origin is given by \eqref{lim.xa}. The proof of Theorem \ref{thm.1.2} relies heavily on the blow-up analysis on $|u_a|$  and the refined energy estimates for $e(a)$ as well. 
Moreover, from the proof of Theorem \ref{thm.1.2}, one can conclude that $e(a)$ satisfies
\begin{equation}\label{lim:ea}
	\frac{e(a)}{(a^\ast-a)^{\frac{1}{2}}}\to\frac{2\lambda^2}{a^\ast}  \text{ as } a  \nearrow a^\ast.
\end{equation}
 We finally remark that Theorem \ref{thm.1.2} also holds for the case where $0<\Lambda<1$ and $\Omega\le2\sqrt{\Lambda}$.

This paper is organized as follows. 
The main purpose of Section \ref{sec.exist} is to prove Theorem \ref{thm.exist} on the existence and non-existence of minimizers for problem \eqref{eq:1.1}. 
Section \ref{sec.MC} is concerned with proving Theorem \ref{thm.1.2} on the limit behavior of minimizers for $e(a)$ as $a\nearrow a^*$.

\section{Existence and non-existence of minimizers}\label{sec.exist}
In this section, we shall complete the proof of Theorem \ref{thm.exist} on the existence and nonexistence of minimizers for $e(a)$. 
Note from \eqref{eq:1.3} and \eqref{ineq:GN.R2} (see also \cite [Lemma 8.1.2] {CT2003}) that,   $w(x)$,  the unique positive solution of  \eqref{eq:1.3},  satisfies
\begin{equation}\label{eq:1.4}
	\int_{\R^2} |\nabla w|^2dx=\int_{\R^2} |w|^2dx=\frac{1}{2} \int_{\R^2} |w|^4dx.
\end{equation}
Moreover, $w(x)$ has the following exponential decay \cite[Proposition 4.1] {GNN}
\begin{equation}\label{eq:1.5}
	w(x), \ |\nabla w(x)|=O(|x|^{-\frac{1}{2}}e^{-|x|})\,\ \text{as} \,\ |x|\rightarrow \infty.
\end{equation}
With the above properties of $w(x)$, we now give the proof of Theorem \ref{thm.exist}. 

\noindent \textbf{Proof of Theorem \ref{thm.exist}.}   
We first prove that $e(a)$ admits at least one minimizer if $0\leq\Omega<\infty$ and $0< a<a^\ast$.  
For any $u\in H^1_0(\mathcal{D},\C)$ satisfying $\|u\|_2^2=1$,  using the GN inequality \eqref{eq:1.6},  one can deduce from \eqref{eq:1.2} and \eqref{ineq: dia.pri} that 
\begin{equation}\label{low:Ea}
\begin{split}
			E_a(u)
		\geq&\int_{\mathcal{D}} |\nabla u|^2dx-\frac{a}{a^*}\int_{\mathcal{D}} |\nabla |u||^2dx-\Omega \int_{\mathcal{D}} x^\bot\cdot\big(iu,\nabla u \big)dx
		\\
		\geq&\frac{a^\ast-a}{a^\ast}\int_{\mathcal{D}} |\nabla u|^2dx-\Omega \int_{\mathcal{D}} x^\bot\cdot\big(iu,\nabla u \big)dx. 
\end{split}
\end{equation}
On the other hand, applying the Cauchy inequality, one has 
\begin{equation}\label{eq:2.2}
	\begin{split}
		\Omega\int_{\mathcal{D}} x^\bot\cdot(iu,\nabla u)dx
		\leq&  \frac{a^\ast-a}{2a^\ast} \int_{\mathcal{D}} |\nabla u|^2dx+\frac{a^\ast\Omega^2}{2(a^\ast-a)}  \int_{\mathcal{D}} |x|^2|u|^2dx\\
		\leq& \frac{a^\ast-a}{2a^\ast} \int_{\mathcal{D}} |\nabla u|^2dx+C,
	\end{split}
\end{equation}
where $C>0$ is a constant depending on $a$, $\Omega$ and $diam(\mathcal{D})$.  
Substituting \eqref{eq:2.2} into \eqref{low:Ea} yields that 
\begin{equation}\label{eq:2.3.0}
	E_a(u)\geq \frac{a^\ast-a}{2a^\ast} \int_{\mathcal{D}} |\nabla u|^2dx-C\geq -C, 
\end{equation} 
which implies that $E_a(u)$ is uniformly bounded from below, and $e(a)$ admits at least one minimizing sequence.

Let $\{u_n\}\subset H^1_0(\mathcal{D},\C) $ be a minimizing sequence of $e(a)$ with $a\in(0,a^*)$, $i.e.$, 
\[
\text{$\|u_n\|^2_2=1$ and $\lim\limits_{n\rightarrow\infty}E_a(u_n)=e(a)$.}
\] 
It then follows from \eqref{eq:2.3.0} that
\begin{equation}\label{eq:2.3}
e(a)\geq \frac{a^\ast-a}{2a^\ast} \int_{\mathcal{D}} |\nabla u_n|^2dx-C \text{ as } n\to\infty,
\end{equation}
which indicates that $\int_{\mathcal{D}} |\nabla u_n|^2dx$ is bounded uniformly in $n$. 
In view of the fact that the embedding $H^1_0(\mathcal{D},\C)\hookrightarrow L^q(\mathcal{D},\C)$ is compact for any $q\in[2,\infty)$ (cf.  \cite[Theorem 1.3.4]{CT2003}),  passing to a subsequence if necessary, there exists some $u_0\in  H^1_0(\mathcal{D},\C)$ such that 
\begin{equation}\label{eq:2.4}
\text{$u_n\overset{n}{\rightharpoonup} u_0$ weakly in $H^1_0(\mathcal{D},\C)$, 
	and 
$u_n\overset{n}{\rightarrow} u_0$ strongly in $L^q(\mathcal{D},\C)$.  }
\end{equation}

Furthermore, one can deduce from \eqref{exp:Eau.2} and \eqref{eq:2.4}  that $E_a(u_n)$ is weak lower semi-continuous on $u_n$, which implies that
\begin{equation}\label{wlsc.Eau}
	e(a)\leq E_a(u_0)\leq \lim\limits_{n\to\infty}E_a(u_n)=e(a). 
\end{equation}
It then follows that $\int_{\mathcal{D}} |u_0|^2dx=1$ and $ e(a)=E_a(u_0)$. 
Hence, for any $0\leq\Omega<\infty$ and $0< a<a^\ast$,  there exists at least one minimizer for $e(a)$.

2.
We next show the non-existence of minimizer for $e(a)$ with $a\geq a^\ast$ by taking some suitable test function.  
Define 
\begin{equation}\label{eq:2.5}
	u_\tau(x):=\frac{A_\tau \tau}{\|w\|_2}\varphi\Big(\frac{x-x_\tau}{R_\tau}\Big)w(\tau(x-x_\tau))e^{i\frac{\Omega}{2}x \cdot x_\tau^\bot},\  x\in \mathcal{D},
\end{equation}
where $\tau >0$, 
$R_\tau>0$, $x_\tau\in\mathcal{D}$ is some points sequence to be determined later, 
$x_\tau^\bot=(-{x_\tau}_{2}, {x_\tau}_{1})$, $A_\tau>0$ is chosen such that $\int_{\mathcal{D}} |u_\tau|^2dx=1$, $\varphi\in C_0^\infty(\R^2)$ is a non-negative smooth cut-off function satisfying 
$\varphi(x)=1$ as  $|x|\leq1$, $\varphi(x)\in(0,1)$ as  $|x|\in(1,2)$ and $\varphi(x)=0$ as  $|x|\geq2$. 

For the convenience of following calculations, here we introduce some new notations and add some necessary conditions on $x_\tau$. 
Denote 
\[
\text{$\mathcal{D}_\tau:=\{\tau(x-x_\tau): x \in\mathcal{D}\}$, and $R_\tau:=\frac{M\ln\tau}{\tau}$, }
\]
where $M>0$ is a constant large enough. 
Take $x_0$ as a minimum point of  $V_\Omega(x)$ in $\bar{\mathcal{D}}$, $i.e.$, 
\begin{equation}\label{inf: x0.V}
	V_\Omega(x_0)=\inf\limits_{x\in\bar{\mathcal{D}}}V_\Omega(x),
\end{equation}
 where $V_\Omega(x)$ is defined by \eqref{def:V_Omega}. 
If  $V_\Omega(x)$ can achieve its minimum at some inner points of $\mathcal{D}$,  set $x_0\in \mathcal{D}$ and $x_\tau\equiv x_0$. It is obvious that $B_{2R_\tau}(x_\tau)\subset \mathcal{D}$ as $\tau$ is large enough.  
If $V_\Omega(x)$ only admits minimum points on the boundary of $\mathcal{D}$, we suppose $x_0\in \partial\mathcal{D}$ and set $x_\tau=x_0-2R_\tau \vec{n}$, where  $\vec{n}$ is the outer unit normal vector of $\partial\mathcal{D}$ at $x_0$. Since $\mathcal{D}$ satisfies the interior ball condition, one can further check that $x_0\in\partial B_{2R_\tau}(x_\tau)$,  $B_{2R_\tau}(x_\tau)\subset \mathcal{D}$, $R_\tau\to 0$ and $x_\tau\to x_0$ as $\tau\to\infty$.

Applying \eqref{eq:1.4} and \eqref{eq:1.5}, 
some calculations yield that 
\begin{equation}\label{eq:2.6}
A^2_\tau=1+o(\tau^{-M}) \text{ as }\tau\rightarrow\ \infty, 
\end{equation}
\begin{equation}\label{eq:2.7}
  \begin{split}
    &\int_{\mathcal{D}} \Big|\Big(\nabla-i\frac{\Omega}{2}x^\bot\Big)u_\tau\Big|^2dx\\
   =&\frac{A_\tau^2 \tau^2}{\|w\|_2^2} \int_{\mathcal{D}}\Big|\frac{1}{R_\tau}\nabla\varphi\Big(\frac{x-x_\tau}{R_\tau}\Big)w(\tau(x-x_\tau))+\tau\varphi\Big(\frac{x-x_\tau}{R_\tau}\Big)\nabla w(\tau(x-x_\tau))\\
    &-i\frac{\Omega}{2}(x-x_\tau)^\bot \varphi\Big(\frac{x-x_\tau}{R_\tau}\Big)w(\tau(x-x_\tau))\Big|^2dx\\
    =&\frac{A_\tau^2\tau^2}{\|w\|_2^2} \int_{\mathcal{D}_\tau} \Big[\varphi^2\Big(\frac{x}{\tau R_\tau}\Big)|\nabla w|^2 +2\frac{1}{\tau R_\tau}\varphi\Big(\frac{x}{\tau R_\tau}\Big)w\nabla\varphi\Big(\frac{x}{\tau R_\tau}\Big)\cdot\nabla w \\
    &+\frac{1}{\tau^2R_\tau^2}  \Big|\nabla \varphi\Big(\frac{x}{\tau R_\tau}\Big)\Big|^2w^2\Big]dx
    +\frac{A_\tau^2}{\|w\|_2^2}  \frac{\Omega^2}{4\tau^2}\int_{\mathcal{D}_\tau}|x|^2\varphi^2\Big(\frac{x}{\tau R_\tau}\Big)w^2 dx\\
    =&\tau^2+\frac{\Omega^2}{4a^*\tau^2} \int_{\R^2} |x|^2w^2dx+o(\tau^{-M}), \text{ as } \tau \rightarrow \infty,
  \end{split}
\end{equation}
\begin{equation}\label{eq:2.9}
	\begin{split}
		\int_{\mathcal{D}} V_\Omega(x)|u_\tau|^2dx
		=&\frac{A^2_\tau }{\|w\|^2_2}\int_{\mathcal{D}_\tau} V_\Omega\Big(\frac{x}{\tau }+x_\tau\Big)  \varphi^2\Big(\frac{x}{\tau R_\tau}\Big) w^2dx\\
		=&V_\Omega(x_0)+o(\tau^{-M}), \text{ as } \tau \rightarrow \infty,
	\end{split}
\end{equation}
and
\begin{equation}\label{eq:2.8}
  \begin{split}
   \frac{a}{2} \int_{\mathcal{D}} |u_\tau|^4dx
  =&\frac{a}{2} \frac{A_\tau^4 \tau^2}{\|w\|_2^4} \int_{\mathcal{D}_\tau} \varphi^4\Big(\frac{x}{\tau R_\tau}\Big)w^4dx\\
    =&\frac{a}{a^*}\tau^2+o(\tau^{-M}), \text{ as } \tau \rightarrow \infty.
  \end{split}
\end{equation}
Using the above estimates, one can deduce from \eqref{exp:Eau.2} that
\begin{equation}\label{eq:2.10}
  e(a)\leq E_a(u_\tau)=\frac{a^\ast-a}{a^\ast}\tau^2+V_\Omega(x_0)+O(\tau^{-2}),  \text{ as } \tau \rightarrow \infty. 
\end{equation}
As $\Omega\geq0$ and $a>a^\ast$, one can directly obtain that $e(a)$ is unbounded from below, and admits no minimizers.

Finally, we shall prove the nonexistence of minimizers for $e(a^*)$ by contradiction. 
It directly follows from \eqref{inf: x0.V} and \eqref{eq:2.10} that
\begin{equation}\label{eq:2.11}
e(a^\ast)\leq V_\Omega(x_0)=\inf\limits_{x\in\bar{\mathcal{D}}} V_\Omega(x).
\end{equation}
 On the other hand, by the GN inequality \eqref{eq:1.6} and the diamagnetic inequality \eqref{eq:1.7}, one can deduce from \eqref{exp:Eau.2} that, for any $u\in H^1_0(\mathcal{D},\C)$ with $\|u\|_2^2=1$, 
\begin{equation}\label{eq:2.12}
  E_{a^*}(u)\geq
  \int_{\mathcal{D}} V_\Omega(x)|u|^2dx\geq\inf\limits_{x\in\bar{\mathcal{D}}} V_\Omega(x).
\end{equation}
Combining \eqref{eq:2.11} and \eqref{eq:2.12} then yields that
\begin{equation}\label{eq:2.13}
   e(a^\ast)=\inf\limits_{x\in\bar{\mathcal{D}}} V_\Omega(x).
\end{equation}
Suppose that there exists a minimizer $u^\ast\not\equiv0$ for $e(a^*)$. Similar to   \eqref{eq:2.12}, one has 
\begin{equation}\label{eq:2.14}
  e(a^\ast)=E_{a^\ast}(u^\ast)\geq\int_{\mathcal{D}} V_\Omega(x)|u^\ast|^2dx>\inf\limits_{x\in\mathcal{D}} V_\Omega(x),
\end{equation}
which contradicts to \eqref{eq:2.13}. This proves the nonexistence of minimizers  for $e(a^\ast)$.  
Hence, the proof of Theorem \ref{thm.exist} is  complete. 
\qed
\vskip 0.05truein

\section{Mass concentration as $a\nearrow a^\ast$}\label{sec.MC}
This section is devoted to completing the proof of Theorem \ref{thm.1.2} by employing the blow-up analysis and energy methods. 
We first address the following energy estimates of $e(a)$ as $a\nearrow a^\ast$.

\begin{lem}\label{lem:engery.est}
Suppose the origin is a inner point of $\mathcal{D}$, and 
$V(x)$ satisfies \eqref{eq:1.8} with $\Lambda>1$ and $\Omega\leq2$.  Then we have
\begin{equation}\label{eq:3.1}
0\leq e(a)\leq  \frac{2\lambda^2+o(1)}{a^\ast}(a^\ast-a)^\frac{1}{2},
\  \mbox{  as }  a  \nearrow  a^\ast.
\end{equation}
where   $\lambda>0$ is defined by \eqref{eq:1.10}.
\end{lem}

\noindent \textbf{Proof.} We start with the lower bound of $e(a)$. 
Using the GN inequality \eqref{eq:1.6} and the diamagnetic inequality \eqref{ineq: dia.pri}, one can  deduce from \eqref{exp:Eau.2}  that, for any $u\in H^1_0(\mathcal{D},\C)$ with ${\|u\|^2_2}=1$,
\begin{equation}
\begin{split}
E_a(u)
&\geq \int_{\mathcal{D}} |\nabla|u||^2dx-\frac{a}{2}\int_{\mathcal{D}} |u|^4dx+\int_{\mathcal{D}} V_\Omega(x)|u|^2dx\\
&\geq \frac{a^\ast-a}{2}\int_{\mathcal{D}} |u|^4dx+\int_{\mathcal{D}} V_\Omega(x)|u|^2dx\geq 0,\\
\end{split}
\end{equation}
which gives the lower bound of \eqref{eq:3.1}.

Next, we shall derive the upper bound in \eqref{eq:3.1} by 
taking a  test function as  \eqref{eq:2.5} with $x_\tau\equiv 0$.  
Similar to \eqref{eq:2.7}-\eqref{eq:2.8}, one has
\begin{equation*}
  \frac{a}{2} \int_{\mathcal{D}} |u_\tau|^4dx
=\frac{a}{a^*}\tau^2+o(\tau^{-M}), \text{ as } \tau \rightarrow \infty,
\end{equation*}
\begin{equation*}
   \int_{\mathcal{D}} \Big|\Big(\nabla-i\frac{\Omega}{2}x^\bot\Big)u_\tau\Big|^2dx
    =\tau^2+\frac{\Omega^2}{4a^*\tau^2} \int_{\R^2} |x|^2w^2dx+o(\tau^{-M}), \text{ as } \tau \rightarrow \infty,
\end{equation*}
and
\begin{equation*}
	\begin{split}
		&\int_{\mathcal{D}} V_\Omega(x)|u_\tau|^2dx\\
		=&\frac{A^2_\tau }{\|w\|^2_2}\int_{\mathcal{D}_\tau}  V_\Omega\Big(\frac{x}{\tau }\Big)   \varphi^2\Big(\frac{x}{\tau R_\tau}\Big) w^2dx\\
		=& \frac{1 }{a^*\tau^2}\int_{\R^2} V_\Omega(x)w^2dx+o(\tau^{-M}), 
		\text{ as } \tau \rightarrow \infty. 
	\end{split}
\end{equation*}
With the above estimates,   take $\tau=\lambda(a^\ast-a)^{-\frac{1}{4}}$ and one has 
\begin{equation}\label{eq:3.3}
  \begin{split}
    e(a)
    \leq&\frac{\tau^2}{a^\ast}(a^\ast-a)+\frac{1}{a^*\tau^2}\int_{\R^2} \Big[V_\Omega(x)+\frac{\Omega^2}{4}|x|^2\Big]w^2dx+o(\tau^{-M})\\
    =&\frac{\tau^2}{a^\ast}(a^\ast-a)+\frac{\lambda^4}{a^\ast \tau^2}+o(\tau^{-M})\\
    =&\frac{2\lambda^2+o(1)}{a^\ast}(a^\ast-a)^\frac{1}{2}, 
    \text{ as } \tau \rightarrow \infty.  
  \end{split}
\end{equation}
We thus obtain the upper bound in \eqref{eq:3.1}, and complete the proof of this lemma. 
\qed

\begin{lem}
	Suppose all the assumptions in Theorem \ref{thm.1.2} hold. Let $u_a$ be a complex-valued minimizer of $e(a)$. Define
	\begin{equation}\label{eq:3.6}
		\epsilon_a:=\Big(\int_{\mathcal{D}} |\nabla u_a|^2dx\Big)^{-\frac{1}{2}}>0.
	\end{equation}
Then we have  
	\begin{enumerate}
		\item  $\epsilon_a\rightarrow0$ and $\mu_a\epsilon_a^2\rightarrow -1$ as $a\nearrow a^\ast$.
		\item  $|u_a|$ admits at least one global  maximal point in $\mathcal{D}$. 
	\end{enumerate}
\end{lem}
\noindent \textbf{Proof.} 
1. On the contrary, suppose $\int_{\mathcal{D}} |\nabla u_a|^2dx\leq C$ for some finite constant $C$. 
Since $\int_{\mathcal{D}} |u_a|^2dx=1$, one can obtain that $\{u_a\}$  is bounded uniformly in $H^1_0(\mathcal{D},\C)$ as $a\nearrow a^\ast$. 
In view of the fact that the embedding $H^1_0(\mathcal{D},\C)\hookrightarrow L^q(\mathcal{D},\C)$ is compact for any $q\in[2,\infty)$ (cf.  \cite[Theorem 1.3.4]{CT2003}),  passing to a sub-sequence (still denoted by $\{u_a\}$) if necessary, there holds that,   
\begin{equation*}
	\begin{split}
		u_a\rightharpoonup u_0\  \mbox{  weakly in } \  H^1_0(\mathcal{D},\C), \ \ u_a\rightarrow u_0\  \mbox{  strongly in } \  L^q(\mathcal{D},\C)
	\end{split}
\end{equation*}
for some $u_0\in H^1_0(\mathcal{D},\C)$. 

Take $\{u_a\}$ with $a\nearrow a^*$ as a minimizing sequence of $e(a^*)$. 
Then some arguments similar to \eqref{eq:2.4}-\eqref{wlsc.Eau} yields that 
$u_0$ is a minimizer of  $e(a^*)$, which contradicts to the fact that $e(a^*)$ admits no minimizer.  
Hence, one has $\epsilon_a\rightarrow 0$ as $a\nearrow a^\ast$.

We next prove that $\mu_a\epsilon_a^2\rightarrow-1$ as $a\nearrow a^\ast$. 
Since $e(a)\to0$ as $a\nearrow a^*$, one can deduce from \eqref{exp:Eau.2} that 
\begin{equation}\label{lim:V}
\int_{\mathcal{D}}V_\Omega(x)|u_a|^2dx\to0,\,\text{ as } a\nearrow a^\ast,  
\end{equation}
and 
\begin{equation}\label{lim:nabla.4}
	\frac{\int_{\mathcal{D}} |u_a|^4dx}{\int_{\mathcal{D}}\Big |\Big(\nabla-i \frac{\Omega}{2}x^{\bot}\Big)u_a\Big|^2dx}\to \frac{2}{a^*},\,\text{ as } a\nearrow a^\ast. 
\end{equation}
On the other hand, using the Cauchy inequality, for any given $\epsilon>0$, there holds that 
\begin{equation}\label{eq:3.15}
	\begin{split}
		\Big|\Omega \int_{\mathcal{D}} x^\bot(i u_a,\nabla u_a)dx\Big|
		\leq& \epsilon \int_{\mathcal{D}} |\nabla u_a|^2dx+\frac{1}{4\epsilon}\int_{\mathcal{D}} \Omega^2 |x|^2|u_a|^2dx. 
	\end{split}
\end{equation}
Take $\epsilon=\Big( \int_{\mathcal{D}} |\nabla u_a|^2dx\Big)^{-\frac12}$, and one has
\begin{equation}\label{eq:3.15.2}
	\Big|\Omega \int_{\mathcal{D}} x^\bot(i u_a,\nabla u_a)dx\Big|\leq C\Big( \int_{\mathcal{D}} |\nabla u_a|^2dx\Big)^\frac12,  
\end{equation}
where $C>0$ is a constant depending on $a$, $\Omega$ and $diam(\mathcal{D})$.  
Furthermore, one can derive from \eqref{eq:1.7} and \eqref{eq:3.15.2} that 
\begin{equation}\label{lim:nabla.ome}
\begin{split}
	&	\frac{\int_{\mathcal{D}}\Big |\Big(\nabla-i \frac{\Omega}{2}x^{\bot}\Big)u_a\Big|^2dx}{\int_{\mathcal{D}} |\nabla u_a|^2dx}\\
	=&1-\frac{\Omega\int_{\mathcal{D}} x^\bot(i u_a,\nabla u_a)dx}{\int_{\mathcal{D}} |\nabla u_a|^2dx}
	+\frac{\frac{\Omega^2}{4 }\int_{\mathcal{D}}  |x|^2|u_a|^2dx
	}{\int_{\mathcal{D}} |\nabla u_a|^2dx}\\
	=&1+o(1),\, \text{ as } a\nearrow a^\ast.  
\end{split}
\end{equation}
It then follows from \eqref{eq:3.5},  \eqref{eq:3.1}, \eqref{lim:nabla.4} and \eqref{lim:nabla.ome} that
\begin{equation}\label{eq:3.13}
\begin{split}
		\lim\limits_{a\nearrow a^\ast} \mu_a\epsilon_a^2
	=&\lim\limits_{a\nearrow a^\ast} \frac{e(a)-\frac{a}{2}\int_{\mathcal{D}} |u_a|^4dx}{\int_{\mathcal{D}} |\nabla u_a|^2dx} \\
	=&\lim\limits_{a\nearrow a^\ast} \frac{e(a)-\frac{a}{2}\int_{\mathcal{D}} |u_a|^4dx}{\int_{\mathcal{D}}\Big |\Big(\nabla-i \frac{\Omega}{2}x^{\bot}\Big)u_a\Big|^2dx}
	=-1. 
\end{split}
\end{equation}

2. 
Define $v_a(x):=\epsilon_a u_a(\epsilon_ax)$. Then one can check that 
\begin{equation}\label{bd:va}
	\int_{\mathcal{D}_a'} |v_a|^2dx=\int_{\mathcal{D}_a'} |\nabla v_a|^2dx  =1,  
\end{equation}
and $v_a(x)$ satisfies 
\begin{equation}\label{equ:va}
	\begin{cases}
		-\Delta v_a+i  \epsilon_a^2 \,\Omega x^\bot\cdot\nabla v_a\\
		+\Big[\frac{\epsilon_a^4 \Omega^2}{4}|x|^2+\epsilon_a^2V_\Omega(\epsilon_a x)
		-\mu_a\epsilon_a^2-a|v_a|^2\Big]v_a=0, \  \text{  in } \mathcal{D}'_a,\\
		v_a=0 ,                \text{  on }  \partial\mathcal{D}'_a, 
	\end{cases}
\end{equation}
where  $\mathcal{D}'_a:=\{\epsilon_a x: x\in\mathcal{D}\}$. 

Suppose $v _a:=R_{v_a}+iI_{v_a}$, where $R_{v_a}$ denotes the real part of $v_a$ and $I_{v_a}$ denotes the imaginary part of $v_a$. 
Then one can further derive that  
\begin{equation}\label{equ:|ua|}
	\begin{cases}
		-\Delta R_{v_a}-  \epsilon_a^2 \,\Omega x^\bot\cdot\nabla I_{v_a}\\
		+\Big[\frac{\epsilon_a^4 \Omega^2}{4}|x|^2+\epsilon_a^2V_\Omega(\epsilon_a x)
		-\mu_a\epsilon_a^2-a|v_a|^2\Big]R_{v_a}=0, \  \text{  in } \mathcal{D}'_a,\\
		-\Delta I_{v_a}+  \epsilon_a^2 \,\Omega x^\bot\cdot\nabla R_{v_a}\\
		+\Big[\frac{\epsilon_a^4 \Omega^2}{4}|x|^2+\epsilon_a^2V_\Omega(\epsilon_a x)
		-\mu_a\epsilon_a^2-a|v_a|^2\Big]I_{v_a}=0, \  \text{  in } \mathcal{D}'_a. 
	\end{cases}
\end{equation}
Using the $L^p$ theory (see e.g., Theorem 8.8 and Theorem 9.11 in \cite{GT}) and  the  standard Sobolev embedding theorem (cf. \cite[Corollary 7.11]{GT}), one can deduce that $R_{v_a}, I_{v_a}\in C^1( \mathcal{D}'_a)$, which implies that $|v_a|^2\in C^1( \mathcal{D}'_a)$. 
Now we shall prove that $|v_a|$ admits at least one global maximum point in $\mathcal{D}_a$ by contradiction. Otherwise, suppose $|v_a|$ achieve its maximum on the boundary of $\mathcal{D}_a$. 
Since $v_a=0$ on  $\partial\mathcal{D}'_a$, one can derive that $|v_a|\equiv0$ in  $\bar{\mathcal{D}}'_a$, which contradicts to \eqref{bd:va}. One can thus directly conclude that $|u_a|$ admits at least one global  maximal point in $\mathcal{D}$. Hence, we complete the proof of this lemma. 
\qed

\begin{lem}\label{lem:wa.1}
Suppose that all the assumptions in Theorem 1.2 hold. Let $u_a$ be a complex-valued minimizer of $e(a)$. Define
\begin{equation}\label{eq:3.7}
	w_a(x):=\epsilon_a u_a(\epsilon_a x+x_a)e^{i(\theta_a-\frac{\epsilon_a \Omega}{2} x \cdot x_a^\bot)},\ \ \ x\in\mathcal{D}_a,
\end{equation}
where $x_a$ is a global maximal point of $|u_a|$, $\theta_a\in[0,2\pi)$ is a proper constant to be determined, and 
\begin{equation}\label{eq:3.8}
	\mathcal{D}_a:=\{x\in\R^2: \epsilon_a x+x_a\in\mathcal{D}\}. 
\end{equation}
Then we have the following. 
\begin{enumerate}
\item  There exist constants $\eta>0$ and $ R>0$, independent of $a$, such that
\begin{equation}\label{eq:3.9}
  \liminf\limits_{a\nearrow a^\ast} \int_{B_{2R}(0)\cap\mathcal{D}_a} |w_a|^2dx\geq\eta >0.
\end{equation}
\item  For each $|u_a|$ with $a$ close enough to $a^*$, the global maximal point $x_a$ is the unique local maximum point of  $|u_a|$,  and satisfies $\lim\limits_{a\nearrow a^\ast} |x_a|=0$.  Moreover, $\mathcal{D}_a$ satisfies
\begin{equation}\label{eq:3.10}
	\mathcal{D}_{a^\ast}=\R^2. 
\end{equation}
and $w_a$  satisfies
\begin{equation}\label{eq:3.11}
  w_a(x)\rightarrow\frac{w(x)}{\sqrt{a^\ast}} \text{ strongly in }  H^1(\R^2,\C)  \text{ as }  a\nearrow a^\ast,
\end{equation}
 where $w$ is the unique positive radially symmetric solution of \eqref{eq:1.3}.
\end{enumerate}
\end{lem}
\noindent \textbf{Proof.} 1. Set $w_a(x):=e^{i\theta_a}\overline{w}_a(x)$, where $\overline{w}_a(x):=\epsilon_a u_a(\epsilon_ax+x_a)e^{-i\frac{\epsilon_a\Omega}{2}x \cdot x_a^\bot}$, and  $\theta_a\in[0,2\pi)$ is chosen properly such that
\begin{equation}\label{eq:3.17}
  \Big\|w_a-\frac{w}{\sqrt{a^\ast}}\Big\|_{L^2(\mathcal{D}_a)}=\min\limits_{\theta\in[0,2\pi)}  \Big\|e^{i\theta} \overline{w}_a-\frac{w}{\sqrt{a^\ast}}\Big\|_{L^2(\mathcal{D}_a)}.
\end{equation}
Then some calculations yield the following orthogonality conditions for $w_a$: 
\begin{equation}\label{eq:3.18}
\int_{\mathcal{D}_a} w(x)I_a(x) dx=0,
\end{equation}
where $I_a(x)$ is the corresponding imaginary part of $w_a(x)$. 

On the other hand, similar to \eqref{lim:V}, \eqref{lim:nabla.4}, \eqref{lim:nabla.ome} and \eqref{bd:va}, one has
\begin{equation}\label{eq:3.19}
  \int_{\mathcal{D}_a} |w_a|^2dx =1,
\end{equation}
\begin{equation}\label{eq:3.20}
	\begin{split}
		\int_{\mathcal{D}_a} |\nabla w_a|^2dx              
		=\epsilon_a^2\int_{\mathcal{D}} \Big|\Big(\nabla -i\frac{ \Omega}{2} x_a^\bot \Big)u_a\Big|^2dx
		\to1, \text{ as }  a  \nearrow  a^\ast,  
	\end{split}
\end{equation}
\begin{equation}\label{lim: wa4}
	\int_{\mathcal{D}_a} |w_a|^4dx \to\frac{2}{a^*}, \text{ as }  a  \nearrow  a^\ast,  
\end{equation}
and
\begin{equation}\label{lim:Vwa}
\int_{\mathcal{D}_a} V_\Omega(\epsilon_a x+x_a)|w_a|^2dx  \to0, \text{ as }  a  \nearrow  a^\ast.   
\end{equation} 
Therefore, $\{w_a\}$ is bounded uniformly in $H_0^1(\mathcal{D}_a,\C)$. 
Using  the diamagnetic inequality \eqref{eq:1.7},  one can check that $|w_a|$ is also bounded uniformly in $H_0^1(\mathcal{D}_a,\R)$. 

Following from \eqref{eq:3.4} and \eqref{eq:3.7}, it can be deduced that $w_a$ is a complex-valued solution of the following Euler-Lagrange equation
\begin{equation}\label{eq:3.21}
\begin{cases}
	 -\Delta\,w_a+i \,\epsilon_a^2 \,\Omega (x^\bot,\nabla w_a)\\
 +\Big[\frac{\epsilon_a^4 \Omega^2}{4}|x|^2+\epsilon_a^2V_\Omega(\epsilon_a x+x_a)
  -\mu_a\epsilon_a^2-a|w_a|^2\Big]w_a=0 \  \mbox{  in } \  \mathcal{D}_a,\\
w_a=0                     \ \  \mbox{  on } \ \ \partial\mathcal{D}_a.
\end{cases}
\end{equation}
Define $W_a:=|w_a|^2\geq0$ in $\mathcal{D}_a$. It then follows that 
\begin{equation}\label{eq:3.22}
\begin{split}
  &-\frac{1}{2}\Delta\,W_a+|\nabla w_a|^2- \epsilon_a^2 \Omega x^\bot (i w_a \cdot \nabla w_a)+ \\
  &\Big[\frac{\epsilon_a^4 \Omega^2}{4}|x|^2+\epsilon_a^2V_\Omega(\epsilon_a x+x_a)-\mu_a\epsilon_a^2-aW_a\Big]W_a=0\  \mbox{  in } \  \mathcal{D}_a.
\end{split}
\end{equation}
From \eqref{eq:1.7}, one has
\begin{equation}\label{eq:3.23}
  |\nabla w_a|^2- \epsilon_a^2 \Omega x^\bot (i w_a \cdot \nabla w_a)+\frac{\epsilon_a^4 \Omega^2}{4}|x|^2 W_a \geq 0\  \mbox{  in } \  \mathcal{D}_a,
\end{equation}
which implies from \eqref{eq:3.22} that
\begin{equation}\label{eq:3.24}
-\frac{1}{2}\Delta W_a-\mu_a\epsilon_a^2 W_a-aW_a^2\leq0\  \mbox{  in } \  \mathcal{D}_a.
\end{equation}
Since 0 is a maximal point of $W_a$, one has $-\Delta W_a(0)\geq0$ for all $a<a^\ast$. 
Moreover, due to $\mu_a \epsilon_a^2\rightarrow-1$ as $a\nearrow a^*$, one can deduce from \eqref{eq:3.24} that 
\begin{equation}\label{eq:3.25}
  W_a(0)\geq\frac{1}{2a}> \frac{1}{2a^\ast}>0.
\end{equation}

Here, we claim that  there exist some constant $R>0$  such that 
\begin{equation}\label{est.lWa.eta}
	\frac{1}{2a^\ast}\leq  \max\limits_{B_R(0)} W_a(x) \leq  C\int_{B_{2R(0)}\cap \mathcal{D}_a} |W_a|^2dx.  
\end{equation}
If there exists a constant $R>0$ such that $B_{2R}(0)\subseteq \mathcal{D}_a$ as $a$ is close enough to $a^\ast$,   
 using the De-Giorgi-Nash-Moser theory \cite[Theorem 4.1]{HL}, one can deduce from \eqref{eq:3.24} that 
\begin{equation}\label{eq:3.26}
  \frac{1}{2a^\ast}< W_a(0)=\max\limits_{B_R(0)} W_a(x) \leq C\int_{B_{2R}(0)} |W_a|^2dx=C\int_{B_{2R(0)}\cap \mathcal{D}_a} |W_a|^2dx.
\end{equation}
On the other hand,  if for any constant $R>0$, $B_{2R}(0)\nsubseteqq \mathcal{D}_a$ as $a\nearrow a^\ast$,  set
\[
\widetilde{W_a}:=
\begin{cases}
	W_a  &  x\in \mathcal{D}_a \\
	0  &  x\in \R^2\backslash\mathcal{D}_a, 
\end{cases}
\]
and then $ \widetilde{W_a}$ satisfies
\begin{equation}\label{eq:3.27}
  -\frac{1}{2}\Delta\widetilde{W_a}-\mu_a\epsilon_a^2 \widetilde{W_a}-a\widetilde{W_a}^2\leq0\,  \mbox{  in } \,  \R^2.
\end{equation}
Thus for any constant $R>0$, there holds that 
\[
\sup\limits_{B_R(0)} \widetilde{W_a}=\widetilde{W_a}(0)=W_a(0)> \frac{1}{2a^\ast}. 
\]
Similar to \eqref{eq:3.26}, one then has 
\begin{equation}\label{eq:3.28}
  \frac{1}{2a^\ast}\leq \sup\limits_{B_R(0)}\widetilde{W_a} \leq C\int_{B_{2R}(0)} |\widetilde{W_a}|^2dx=C\int_{B_{2R}(0)\cap \mathcal{D}_a} |W_a|^2dx.
\end{equation}
Hence, \eqref{est.lWa.eta} follows from \eqref{eq:3.26} and \eqref{eq:3.28}.
Furthermore, one can deduce from  \eqref{est.lWa.eta} that, there exists some constant $\eta>0$ such that   
\begin{equation}\label{eq:3.30}
  \int_{B_{2R}(0)\cap \mathcal{D}_a} |w_a|^2dx \geq \int_{B_{2R}(0)\cap \mathcal{D}_a} \frac{W_a^2}{\sup\limits_{x\in \mathcal{D}_a} W_a}dx \geq \eta >0.
\end{equation}
\eqref{eq:3.9} is thus proved.

2. Now we shall prove that $\lim\limits_{a\nearrow a^\ast} |x_a|=0$, where $x_a$ is a global maximal point of $|u_a|$.
By contradiction, suppose  there exists $c_0>0$ such that $|x_a|\geq c_0$ uniformly in $a$. One can derive from  \eqref{eq:3.9} and \eqref{lim:Vwa} that
\begin{equation}\label{eq:3.31}
  \begin{split}
    0=
    &\lim\limits_{a\nearrow a^\ast} \int_{\mathcal{D}_a} V_\Omega(\epsilon_a x+x_a)|w_a|^2dx \\
    \geq& \int_{B_{\frac{1}{\sqrt{\epsilon_a}}(0)}\cap\mathcal{D}_a} \liminf\limits_{a\nearrow a^\ast} \Big(1-\frac{\Omega^2}{4}\Big)|\epsilon_a x+x_a|^2|w_a    |^2dx \\
    \geq & C\eta>0,
  \end{split}
\end{equation}
which is a contradiction. 

Since $\epsilon_a\rightarrow0$ and $|x_a|\rightarrow0\in\mathcal{D}$ as $a\nearrow a^\ast$, one can observe from \eqref{eq:3.8} that $\mathcal{D}_{a^\ast}=\R^2$.
Hence, one can extend the definition domain of $w_a(x)$ to $\R^2$ by setting $w_a\equiv0$ in $\R^2\setminus\mathcal{D}_a$. In the following proof, we shall always suppose $w_a\in H^1(\R^2,\C)$ in the sense that $w_a\equiv0$ in $\R^2\setminus\mathcal{D}_a$. 
Since $|w_a|$ is bounded uniformly in $H^1(\R^2,\R)$, passing to a subsequence if necessary, there holds that  
\[
\text{$|w_a|\rightharpoonup w_0$ weakly in $ H^1(\R^2,\R)$ as $a\nearrow a^\ast$ for some $w_0\in H^1(\R^2,\R)$. }
\]

Next, we shall certify that
\begin{equation}
	\text{$|w_a|\to w_0$ strongly in $ H^1(\R^2,\R)$ as $a\nearrow a^\ast$. }
\end{equation}
Note from \eqref{eq:3.9} that $w_0\not\equiv 0$ in $\R^2$. By the weak convergence, one may assume that $|w_a|\rightarrow w_0$ almost everywhere in $\R^2$ as $a\nearrow a^\ast$.
Applying the Br$\acute{e}$zis-Lieb Lemma (cf. \cite[Lemma 1.2.3]{CPY}) gives that
\begin{equation}\label{eq:3.33}
 \Big\|\nabla |w_a| \Big\|^2_2=\|\nabla w_0\|^2_2+\Big\|\nabla |w_a|-\nabla w_0\Big\|^2_2+o(1), \text{ as $a\nearrow a^*$}, 
\end{equation}
and
\begin{equation}\label{eq:3.32}
	\begin{split}
		\|w_a\|^q_q=\|w_0\|^q_q+\Big\||w_a|-w_0\Big\|^q_q+o(1),  \text{ as $a\nearrow a^*$,}
	\end{split}
\end{equation}
where  $q\in[2,+\infty)$. 
 It then follows that 
 \begin{equation}\label{eq:3.36}
 	\begin{split}
 		 &\lim\limits_{a\nearrow a^\ast} \Big( \int_{\mathcal{D}_a} \big| \nabla |w_a| \big|^2dx-\frac{a^\ast}{2} \int_{\mathcal{D}_a} |w_a|^4dx \Big) \\
 		=&\int_{\R^2} |\nabla w_0|^2dx-\frac{a^\ast}{2} \int_{\R^2} |w_0|^4dx\\
 		&+\lim\limits_{a\nearrow a^\ast} \Big[ \int_{\mathcal{D}_a} \big|\nabla ( |w_a|- w_0)\big|^2dx-\frac{a^\ast}{2} \int_{\mathcal{D}_a} \big| |w_a|-w_0 \big|^4dx \Big] 
 	\end{split}
 \end{equation}
 
On the other hand,  using the diamagnetic inequality \eqref{ineq: dia.pri}, GN inequality \eqref{eq:1.6},    \eqref{eq:3.20}-\eqref{lim: wa4}, one can derive that 
\begin{equation}\label{eq:3.35}
\begin{split}
  0\leq&\lim\limits_{a\nearrow a^\ast}  \Big[\int_{\mathcal{D}_a} \big|\nabla |w_a|\big|^2dx-\frac{a^\ast}{2} \int_{\mathcal{D}_a} |w_a|^4dx \Big] \\
 \leq & \lim\limits_{a\nearrow a^\ast} \Big[\int_{\mathcal{D}_a} |\nabla w_a|^2dx-\frac{a^\ast}{2} \int_{\mathcal{D}_a} |w_a|^4dx\Big]=0. \\
  \end{split}
\end{equation}
Since $\|w_0\|_2^2\leq \lim\limits_{a\nearrow a^*}\|w_a\|_2^2=1$ and $\| |w_a|-w_0\|_2^2\leq 1$,  using the GN inequality \eqref{eq:1.6} again, one can deduce from \eqref{eq:3.36} and \eqref{eq:3.35} that 
\begin{equation}\label{id:w0}
	\int_{\R^2} |\nabla w_0|^2dx=\frac{a^\ast}{2} \int_{\R^2} |w_0|^4dx, 
\end{equation}
which implies from  \eqref{ineq:GN.R2} that  $\|w_0\|_2^2\geq1$. 
Hence, one has $\|w_0\|_2^2=1$ and $|w_a|\rightarrow w_0$ strongly in $L^2(\R^2)$ as $a\nearrow a^\ast$. 
Due to the uniform boundedness of $|w_a|$ in $H^1(\R^2)$, using the interpolation inequality, one has
\begin{equation}\label{lim:|wa|.Lp}
	\text{ $|w_a|\rightarrow w_0$ strongly in $L^q(\R^2 )$ as $a\nearrow a^\ast$, where $ q\in [2,\infty)$. }
\end{equation}
It then follows from \eqref{lim: wa4} and \eqref{eq:3.32} that 
\[
\int_{\R^2} |w_0|^4dx=\lim\limits_{a\nearrow a^\ast}\int_{\R^2} |w_a|^4dx=\frac{2}{a^*}, 
\]
which implies from \eqref{id:w0} that $\int_{\R^2} |\nabla w_0|^2dx =1$. 
Furthermore, one can get from \eqref{eq:3.20} and \eqref{eq:3.33} that 
$\nabla |w_a|\rightarrow \nabla w_0$ strongly in $L^2(\R^2)$ as $a\nearrow a^\ast$. 
Hence, one can finally obtain that 
\[
\text{$|w_a|\rightarrow w_0$ strongly in $H^1(\R^2,\R)$ as $a\nearrow a^\ast$. }
\]

 Note from \eqref{id:w0} that $w_0$ is an optimizer of the GN inequality \eqref{ineq:GN.R2}  in $\R^2$.  Together with the fact $\|w_0\|_2^2=\|\nabla w_0\|_2^2=1$, this indicates that, up to a translation if necessary, $w_0=\frac{1}{\sqrt{a^*}}w(x)$, $i.e.$,  there exists some $y_0\in\R^2$ such that 
 \[
 w_0=\frac{1}{\sqrt{a^*}}w(x+y_0), 
 \]
  where $w$ is the unique positive radially symmetric solution of \eqref{eq:1.3}. 
Thus, one can conclude that,  up to a subsequence if necessary, there holds that 
\begin{equation}\label{eq:3.38}
  |w_a|\rightarrow  \frac{w(x+y_0)}{\sqrt{a^\ast}} \text{ strongly in }  H^1(\R^2,\R)  \text{ as } a\nearrow a^\ast, 
\end{equation}
which indicates that 
\begin{equation}\label{eq:3.39}
	w_a\to\frac{w(x+y_0)}{\sqrt{a^\ast}}e^{i\sigma} \text{ strongly in $H^1(\R^2,\C)$ as $a \nearrow a^*$,} 
\end{equation}
for some $\sigma\in [0,2\pi)$. Using the orthogonality condition \eqref{eq:3.18}, one has $\sigma=0$,

Recall from \eqref{eq:3.21} that $w_a(x)$ satisfies   
\begin{equation}\label{eq:3.47}
	-\Delta w_a=f(w_a)  \text{ in }   H^1(\R^2,\C), 
\end{equation}
where $w_a$ in $\R^2$ is extended from $w_a$ in $\mathcal{D}_a$ by setting $w_a\equiv0$ in $\R^2\setminus \mathcal{D}_a$, and 
\[
f(w_a):=-i \epsilon_a^2 \Omega (x^\bot,\nabla w_a)-\Big[\frac{\epsilon^4_a \Omega^2}{4}|x|^2+V_\Omega(\epsilon_a x+x_a)-\mu_a\epsilon^2_a-a|w_a|^2\Big]w_a. 
\]
Since $w_a$ is bounded uniformly in $H^1(\R^2,\C)$, one can check that $f(w_a)$ is bounded uniformly in $L^2_{loc}(\R^2,\C)$ as $a\nearrow a^*$. 
Using \cite[Theorem 8.8]{GT}, one can then derive that $w_a\in H^2(B_R)$ and  
\begin{equation}\label{eq:3.48}
	\|w_a\|_{H^2(B_R)} \leq C(\|w_a\|_{H^1(B_{R+1})}+\|f(w_a)\|_{L^2(B_{R+1})}), 
\end{equation}
where $R>0$ is a large constant.  
Apply the $L^p$ estimate (cf. \cite[Theorem 9.11]{GT}) and one can deduce that $w_a$ is bounded uniformly in $W^{2,2}(B_R)$. Using the standard Sobolev embedding thus yields that  $w_a$ is bounded uniformly in $C^{1,\alpha}(B_R)$. 
Further, since $V(x)\in C^\infty(\mathcal{D})$ and $\partial\mathcal{D}\in C^1$, by the Schauder's theory, one can deduce that $w_a$ is bounded uniformly in $C^{2,\alpha}(B_R)$ and \eqref{eq:3.38}  holds in $C^2(B_R)$, which implies that 
\begin{equation}\label{lim: wa.w0.C2}
	|w_a|\rightarrow  \frac{w(x+y_0)}{\sqrt{a^\ast}} \text{ strongly in }  C^2_{loc}(\R^2,\R) \text{ as } a\nearrow a^\ast, 
\end{equation} 
because of the arbitrariness of  $R>0$.  
Since  the origin is a global maximal point of  $|w_a|$, so it must also be a global maximal point of $w(x+y_0)$. This  implies from \eqref{lim: wa.w0.C2} that 
\[
y_0=0, 
\]
in view of the fact that the origin is the unique global maximal point of $w$.  
Moreover, since the convergence above is independent of the choice of subsequence, \eqref{eq:3.39} holds true for the whole sequence. 
Hence, \eqref{eq:3.11} is thus proved. 

Finally, we shall prove the uniqueness of local maximum point of $|u_a|$. 
Some arguments similar to 
\eqref{eq:3.22}-\eqref{est.lWa.eta} yield that $|w_a|\to0$ as $|x|\to\infty$  and  $w_a(z_0)\geq\frac{1}{2a^*}>0$ if $a$ is close enough to $a^*$, where $z_0$ is any local maximum point of $|w_a|$. One can thus conclude from \eqref{lim: wa.w0.C2} that all the  local maximal point of  $|w_a|$ must stay in a finite ball as $a\nearrow a^*$. 
Moreover, since  $w(x)$ is radially symmetric about the origin and decreases strictly in $|x|$, by \cite[Lemma 4.2]{NT}, one can deduce that the origin is the unique maximal point of  $|w_a|$ as $a$ is close enough to $a^*$.  This then gives the   uniqueness of local maximum point for $|u_a|$ as $a$ is close enough to $a^*$. 
Hence, we complete the proof of this lemma.
\qed

\begin{lem}\label{lem:wa.2}
Suppose that all the assumptions in Theorem 1.2 hold. 
Let $u_a$ be a complex-valued minimizer of $e(a)$, and $w_a$ is defined by \eqref{eq:3.7}.  
Then we have the following conclusions. 
\begin{enumerate}
\item  There exists a large constant $R\geq0$ such that
\begin{equation}\label{eq:3.40}
  |w_a|\leq Ce^{-\frac{1}{4}|x|} \text{ in }   \R^2\backslash B_R(0).
\end{equation}
\item  $w_a(x)$ satisfies
\begin{equation}\label{eq:3.41}
  w_a(x)\rightarrow \frac{w(x)}{\sqrt{a^\ast}} \text{ uniformly in }  L^\infty(\R^2,\C) \text{ as }  a\nearrow \ a^\ast.
\end{equation}
\item The following estimate holds: 
\begin{equation}\label{eq:3.42}
  \Omega \int_{\mathcal{D}_a} x^\bot(iw_a,\nabla w_a)dx=o(\epsilon _a^2)  \text{ as }  a\nearrow \ a^\ast,
\end{equation}
where $\epsilon_a >0$ is defined by \eqref{eq:3.6}.
\end{enumerate}
\end{lem}
\noindent \textbf{Proof.} 1. Recall from \eqref{eq:3.10} that $\mathcal{D}_{a^\ast}=\R^2$. As in Lemma \ref{lem:wa.1}, we shall extend $w_a(x)$ to $\R^2$ by setting $w_a\equiv0$ in $\R^2\setminus \mathcal{D}_a$. 
Using the H\"{o}lder inequality and the Sobolev inequality, one can deduce from \eqref{eq:3.11} that
\[
W_a\to\frac{w^2}{a^\ast} \text{ strongly in $H^1(\R^2,\R)$ as $a \nearrow a^*$.} 
\]
Then by  the De-Giorgi-Nash-Moser theory \cite[Theorem 4.1]{HL}, for any given $\epsilon'>0$, there exists a sufficiently large $R>0$ such that $W_a(x)<\epsilon'$ in $\R^2\backslash B_R(0)$.  
Moreover, from \eqref{eq:3.24}, one has 
\begin{equation*}%\label{eq:3.45}
	-\Delta W_a+W_a\leq 0\  \mbox{ in } \  \R^2\backslash B_R(0).
\end{equation*}
%%Set $v:=Me^{-\frac{1}{2}|x|}$, where $M>0$ is a constant large enough. And one can check that $-\Delta v+v\geq0$  in $\R^2$.  
%%Furthermore, one can verify that there exist some $R$ and $M$ large enough such that $ W_a\leq Me^{-\frac{1}{2}R} $ on $|x|=R$.  
Thus using  the comparison principle, one can deduce that $|w_a|$ decays  exponentially as \eqref{eq:3.40}.

2. Some arguments similar to \eqref{eq:3.47} and \eqref{eq:3.48} yield that  $w_a$ is uniformly bounded in $ H^2(B_R)$, where  $R>0$ is a large constant. 
Furthermore,  because the embedding $H^2(B_R)\hookrightarrow L^\infty(B_R)$ is compact, one can deduce that,  passing to a subsequence if necessary, \eqref{eq:3.11} holds in $L^\infty(B_R)$. 
This further implies   that \eqref{eq:3.11} holds in $L_{loc}^\infty (\R^2)$, due to that $R>0$ is arbitrary. 
Combining with the exponential decay of  $|w_a|$ in \eqref {eq:3.40}, it then follows that, passing to a subsequence if necessary, \eqref{eq:3.11} holds in $L^\infty (\R^2)$ uniformly.  Moreover, since the convergence is independent of the choice of subsequence, \eqref{eq:3.41} holds true for the whole sequence.

3. Recall from \eqref{eq:1.2} that
\begin{equation}\label{eq:3.51}
  \begin{split}
    e(a)=&\frac{1}{\epsilon^2_a} \int_{\mathcal{D}_a} |\nabla w_a|^2dx-\frac{a}{2 \epsilon^2_a} \int_{\mathcal{D}_a} |w_a|^4dx+\frac{\epsilon^2_a \Omega^2}{4} \int_{\mathcal{D}_a} |x|^2|w_a|^2dx \\
    &+\int_{\mathcal{D}_a} V_\Omega(\epsilon_a x+x_a)|w_a|^2dx-\Omega \int_{\mathcal{D}_a} x^\bot \cdot(i w_a,\nabla w_a)dx \\
    \geq& \frac{1}{\epsilon^2_a} \int_{\mathcal{D}_a} |\nabla w_a|^2dx-\frac{a^\ast}{2 \epsilon^2_a} \int_{\mathcal{D}_a} |w_a|^4dx-\Omega \int_{\mathcal{D}_a} x^\bot (i w_a,\nabla w_a)dx. \\
  \end{split}
\end{equation}
Set $w_a(x)=R_a(x)+iI_a(x)$, where  $R_a(x)$ and $I_a(x)$ denote the real part and imaginary part of $w_a$ respectively. 
It then follows from  \eqref{eq:3.19}, \eqref{eq:3.20} and \eqref{eq:3.41} that 
\[
\int_{\mathcal{D}_a} | w_a|^2dx=\int_{\R^2} | R_a|^2+|  I_a|^2dx=1,
\]
\[
\int_{\mathcal{D}_a} |\nabla w_a|^2dx=\int_{\R^2} |\nabla R_a|^2+|\nabla I_a|^2dx\to1, \text{ as } a\nearrow a^*,  
\]
and 
\begin{equation}\label{lim:Ra.Ia}
	R_a(x)\to \frac{w(x)}{\sqrt{a^\ast}} \text{ and }  I_a(x)\to0 \text{ uniformly  in }  L^\infty(\R^2,\R) \mbox{ as } a\nearrow  a^\ast.
\end{equation}

Moreover, using the GN inequality \eqref{ineq:GN.R2} and \eqref{lim:Ra.Ia},  some calculations yield  
 \begin{equation*}
  \begin{split}
  &\frac{a^\ast}{2}\int_{\mathcal{D}_a} |w_a|^4dx\\
  =&\frac{a^\ast}{2}\int_{\R^2}\big(R_a^4+I_a^4+2R_a^2I_a^2\big)dx\\
  \leq&\int_{\R^2} |\nabla R_a|^2dx \int_{\R^2} R_a^2dx+
  a^\ast\int_{\R^2} R_a^2I_a^2 dx
  +o(1)\int_{\R^2} I_a^2 dx\\
   =&\int_{\R^2} |\nabla R_a|^2dx \Big(1-\int_{\R^2} I_a^2dx\Big)+
   \big((1+o(1)\big)\int_{\R^2}w^2I_a^2dx
   +o(1)\int_{\R^2} I_a^2 dx\\
    =&\int_{\R^2} |\nabla R_a|^2dx-\int_{\R^2}I_a^2dx+ \int_{\R^2}w^2I_a^2dx+o(1)\int_{\R^2}I_a^2dx,  
  \text{ as } a  \nearrow a^\ast, 
  \end{split}
\end{equation*}
and 
\begin{equation}\label{est:omega}
\begin{split}
  &\Omega \int_{\mathcal{D}_a} x^\bot \cdot(i w_a,\nabla w_a)dx\\
  =&\Omega \int_{\R^2} x^\bot \cdot(R_a \nabla I_a-I_a \nabla R_a)dx\\
  =&2\Omega \int_{\R^2} x^\bot\cdot R_a \nabla I_a dx  \\
    =&\frac{2\Omega}{\sqrt{a^*}} \int_{\R^2} x^\bot\cdot w \nabla I_a dx +\frac{2\Omega}{\sqrt{a^*}} \int_{\R^2} x^\bot\cdot (\sqrt{a^*}R_a-w) \nabla I_a dx\\
   \leq&\frac{2\Omega}{\sqrt{a^*}} \int_{\R^2} x^\bot\cdot w \nabla I_a dx +o(1)\|\nabla I_a \|_2 \\
   =&o(1)\|\nabla I_a \|_2,  \text{ as } a  \nearrow  a^\ast,
  \end{split}
\end{equation}
where the last ``$=$'' holds due to  the fact that
\[
\int_{\R^2} x^\bot\cdot w \nabla I_a dx=-\int_{\R^2} x^\bot\cdot  \nabla w I_a dx,
\]
and 
\[
 x^\bot\cdot  \nabla w(x)= x^\bot\cdot  \nabla_x w(|x|)=0, \text{ in } \R^2. 
\]

Combining all the above estimates, one can derive from \eqref{eq:3.51} that
\begin{equation*}\label{eq:3.52}
\begin{split}
		\epsilon_a^2 e(a) 
	\geq& \|I_a\|^2_{H^1(\R^2)} -\int_{\R^2}w^2I_a^2dx+o(1)\| I_a\| _2^2-o(1) \epsilon_a^2\|\nabla I_a \| _2\\
	=& \frac12\|I_a\|^2_{H^1(\R^2)} -\int_{\R^2}w^2I_a^2dx-o(1) \epsilon_a^2\|\nabla I_a \| _2, 
	\text{ as } a  \nearrow a^\ast. 
\end{split}
\end{equation*}
Recall from \eqref{eq:3.18} that $I_a$ is orthogonal to $w$, and one then has $\int_{\R^2}w^2I_a^2dx=0$. It then follows that 
\begin{equation}\label{eq:3.57}
	\epsilon_a^2 e(a)\geq \frac{1}{2}\|I_a\|^2_{H^1(\R^2)}-o(1) \epsilon_a^2\|\nabla I_a \| _2,  \text{ as } a  \nearrow a^\ast. 
\end{equation}

On the other hand, using the GN inequality \eqref{ineq:GN.R2} and the diamagnetic inequality \eqref{eq:1.7}, one can derive from  \eqref {exp:Eau.2}, \eqref{eq:3.1} and \eqref{lim:nabla.ome} that 
\[
\begin{split}
	 \frac{2\lambda^2+o(1)}{a^\ast}(a^\ast-a)^\frac{1}{2}
 \geq	 e(a)
	 \geq& \frac{a^*-a}{a^*} \int_{\mathcal{D}}\Big |\Big(\nabla-i \frac{\Omega}{2}x^{\bot}\Big)u_a\Big|^2dx\\
	  =&\frac{a^*-a}{a^*}(1+o(1))\epsilon_a^{-2},  \text{ as } a  \nearrow a^\ast.  
\end{split}
\]
This indicates that 
\[
  (a^\ast-a)^{\frac{1}{2}}\leq C\epsilon_a^2 \ \text{ and } \ e(a)\leq C\epsilon_a^2,  \text{ as } a  \nearrow a^\ast. 
\]
From \eqref{eq:3.57}, one can then  obtain that 
\[
  \| I_a\|_{H^1(\R^2)} \leq C\epsilon_a^2 \ \text{ and } \  
\|\nabla I_a \| _2\leq C\epsilon_a^2,  \text{ as } a  \nearrow a^\ast.   
\]
Furthermore, from \eqref{est:omega}, one has 
\begin{equation}\label{eq:3.60}
\begin{split}
	\Big|\Omega \int_{\mathcal{D}_a} x^\bot \cdot(i w_a,\nabla w_a)dx\Big|
	=&
\Big|\int_{\R^2} x^\bot\cdot (\sqrt{a^*}R_a-w) \nabla I_a dx\Big|\\
\leq& o(1)  \|\nabla I_a \| _2\leq o(1) \epsilon_a^2. 
\end{split}
\end{equation}
Hence, \eqref{eq:3.42} follows from  \eqref{eq:3.60}. 
We thus complete the proof of this lemma. 
\qed

\noindent \textbf{Proof of Theorem \ref{thm.1.2}}  
With the above lemmas, we shall finally complete the proof  Theorem \ref{thm.1.2} by establishing the refined energy estimates. 
Similar to \eqref{eq:3.51}, from  \eqref{eq:3.19} and \eqref{eq:3.20}, one has 
\begin{equation}\label{eq:3.62}
\begin{split}
  e(a)
  =& \frac{1+o(1)}{\epsilon^2_a} \frac{a^\ast-a}{a^*} 
  +\int_{\mathcal{D}_a} V_\Omega(\epsilon_a x+x_a)|w_a|^2dx\\
  &+\frac{\epsilon^2_a \Omega^2}{4} \int_{\mathcal{D}_a} |x|^2|w_a|^2dx\\
 & -\Omega \int_{\mathcal{D}_a} x^\bot (i w_a,\nabla w_a)dx,  \text{ as } a  \nearrow a^\ast.  
  \end{split}
\end{equation}
Using \eqref{eq:3.41}, one can deduce that
\begin{equation}\label{eq:3.64}
  \frac{\epsilon^2_a \Omega^2}{4} \int_{\mathcal{D}_a} |x|^2|w_a|^2dx=\frac{(1+o(1)) \Omega^2\epsilon^2_a}{4a^\ast}\int_{\R^2}|x|^2w^2dx,  \text{ as } a  \nearrow a^\ast,
\end{equation}
and 
\begin{equation}\label{eq:3.65}
		\int_{\mathcal{D}_a} V_\Omega(\epsilon_a x+x_a)|w_a|^2dx
		=\epsilon_a^2 \int_{\mathcal{D}_a}  V_\Omega\Big(x+\frac{x_a}{\epsilon_a }\Big)|w_a|^2dx. 
\end{equation}

Here, we claim that 
\begin{equation}\label{lim:xa.epsa}
\text{$\frac{|x_a|}{\epsilon_a }\to0$ as $a\nearrow a^*$. }
\end{equation}
We shall first prove that $\frac{|x_a|}{\epsilon_a }$ is bounded uniformly as $a\nearrow a^*$. Otherwise, suppose $\frac{|x_a|}{\epsilon_a }\to+\infty$ as $a\nearrow a^*$, and one can deduce from \eqref{def:V_Omega} and \eqref{eq:3.65} that, there exists some large  constant $M>0$  such that 
\begin{equation}\label{eq:3.65.1}
		\int_{\mathcal{D}_a} V_\Omega(\epsilon_a x+x_a)|w_a|^2dx
	\geq M\epsilon_a^2.  
\end{equation}
It then follows from \eqref{eq:3.42},  \eqref{eq:3.62}-\eqref{eq:3.64} and \eqref {eq:3.65.1}    that
\begin{equation*}
	e(a)
	\geq\frac{a^\ast-a}{a^\ast\epsilon^2_a}+M\epsilon^2_a+o(\epsilon^2_a)
	\geq M'(a^\ast-a)^{\frac{1}{2}},   \text{ as } a  \nearrow a^\ast,  
\end{equation*}
where $M'>0$ is a constant large enough. 
This contradicts to the upper energy estimate in \eqref{eq:3.1}, and  the uniform boundedness  of $\frac{|x_a|}{\epsilon_a }$  as $a\nearrow a^*$ is thus proved. 

Next, we shall complete the proof of  \eqref{lim:xa.epsa}. 
Otherwise, suppose that, passing to a subsequence if necessary,   there exists some  $y_0\neq(0,0)$ such that $\frac{x_a}{\epsilon_a }\to y_0$ as $a\nearrow a^*$.  
Set 
\[
\kappa:= \Big[\int_{\R^2}  V_\Omega(x+y_0)w^2dx+\frac{\Omega^2}{4} \int_{\R^2} |x|^2w^2dx\Big]^\frac{1}{4}. 
\]
And one can check from \eqref{def:V_Omega} that $\kappa>\lambda$, where $\lambda$ is defined by \eqref{eq:1.10}. 
From \eqref{eq:3.41}, \eqref{eq:3.64} and \eqref{eq:3.65}, one can thus  derive  that 
\begin{equation*}  
	\begin{split}
	&\int_{\mathcal{D}_a} V_\Omega(\epsilon_a x+x_a)|w_a|^2dx 
	+ \frac{\epsilon^2_a \Omega^2}{4} \int_{\mathcal{D}_a} |x|^2|w_a|^2dx\\
	=&\big(1+o(1)\big)\frac{\epsilon^2_a}{a^\ast}\int_{\R^2}  V_\Omega(x+y_0)|w|^2dx+\frac{(1+o(1)) \Omega^2\epsilon^2_a}{4a^\ast}\int_{\R^2}|x|^2w^2dx\\
	=&\big(1+o(1)\big)\frac{\epsilon^2_a}{a^\ast}\kappa^4,    
	\text{ as } a  \nearrow a^\ast. 
\end{split}
\end{equation*} 
Substituting this estimates into \eqref{eq:3.62}, it then follows from \eqref{eq:3.42} that 
\begin{equation*}\label{eq:3.kapaa}
	e(a)
	\geq\frac{a^\ast-a}{a^\ast\epsilon^2_a}+\frac{\epsilon^2_a}{a^\ast}\kappa^4+o(\epsilon^2_a)
	\geq\frac{2\kappa^2+o(1)}{a^\ast}(a^\ast-a)^{\frac{1}{2}},   \text{ as } a  \nearrow a^\ast.  
\end{equation*}
This contradicts to the upper energy estimate in \eqref{eq:3.1}, due to the fact $\kappa>\lambda$, and the claim \eqref{lim:xa.epsa} is thus proved.

With the above claim \eqref{lim:xa.epsa}, one can thus  derive from \eqref{eq:3.41} and  \eqref{eq:3.65}  that 
\begin{equation}\label{eq:3.65.3}
	\begin{split}
	&\int_{\mathcal{D}_a} V_\Omega(\epsilon_a x+x_a)|w_a|^2dx\\
	=&\big(1+o(1)\big)\frac{\epsilon^2_a}{a^\ast}\int_{\R^2}  V_\Omega(x)w^2dx\\
	=&\big(1+o(1)\big)\frac{\epsilon^2_a}{a^\ast}\lambda^4-\frac{\epsilon^2_a \Omega^2}{4} \int_{\R^2} |x|^2w^2dx,   \text{ as } a  \nearrow a^\ast,  
	\end{split}
\end{equation}
where $\lambda$ is defined by \eqref{eq:1.10}. 
Following from \eqref{def:V_Omega},   \eqref{eq:3.62}-\eqref{lim:xa.epsa} and \eqref{eq:3.65.3}, one has
\begin{equation}\label{eq:3.66}
		e(a)
		\geq\frac{a^\ast-a}{a^\ast\epsilon^2_a}+\frac{\epsilon^2_a}{a^\ast}\lambda^4+o(\epsilon^2_a)
		\geq\frac{2\lambda^2+o(1)}{a^\ast}(a^\ast-a)^{\frac{1}{2}},   \text{ as } a  \nearrow a^\ast. 
\end{equation}
Together with  the upper energy estimate in \eqref{eq:3.1}, one can directly  obtain that 
\begin{equation}\label{est.ea}
  e(a)=\frac{2\lambda^2+o(1)}{a^\ast}(a^\ast-a)^{\frac{1}{2}},  \text{ as } a  \nearrow a^\ast. 
\end{equation}
which implies all the ``$=$'' in \eqref{eq:3.66} hold,  
and $\epsilon_a $ thus satisfies 
\begin{equation}\label{eq:3.67}
  \epsilon_a =\frac{1+o(1)}{\lambda}(a^\ast-a)^{\frac{1}{4}},  \text{ as } a  \nearrow a^\ast. 
\end{equation}
Moreover, since the convergence above is independent of the choice of subsequence, one can check that \eqref{lim:xa.epsa} and \eqref{eq:3.67} hold true for the whole sequence. 

Finally, set  
\[
\tilde{u}_a:=\sqrt{a^\ast}\frac{(a^\ast-a)^{\frac{1}{4}}}{\lambda} u_a\Big(\frac{(a^\ast-a)^{\frac{1}{4}}}{\lambda} x+x_a\Big)e^{-i(\frac{\Omega}{2\lambda}(a^\ast-a)^{\frac{1}{4}} x \cdot x_a^\bot-\theta_a)}, 
\]
 where  $\theta_a\in[0,2\pi)$ is a proper constant and $\lambda$ is defined by \eqref{eq:1.10}. One can then directly obtain that  \eqref{eq:1.9}  follows from \eqref{eq:3.41} and \eqref{eq:3.67}, \eqref{lim.xa} follows from \eqref{lim:xa.epsa}, and \eqref{lim:ea} follows from \eqref{est.ea}. 
Hence, we complete the proof Theorem \ref{thm.1.2}. 
\qed

\vskip 0.2truein

\noindent {\bf Acknowledgements:} 
This research was  partially supported by 
the National Natural Science Foundation of China (Grant No. 12371113 and 11901223). 
The authors are grateful to Yujin Guo for his fruitful discussions on the present paper.

\end{document}